\documentclass[10pt,fleqn]{amsart}

\usepackage{amsmath,amssymb,latexsym}

\theoremstyle{plain}
\newtheorem*{theorem}{Theorem}

\theoremstyle{remark}
\newtheorem*{ack}{Acknowledgment}

\numberwithin{equation}{section}

\def\mydate{\number\year-\ifnum\month<10{0}\fi\number\month-\ifnum\day<10{0}\fi\number\day}

\newcommand{\dy}{\partial}
\newcommand{\ddt}[1]{\frac{\mathrm{d}{#1}}{\mathrm{d}{t}}}

\newcommand{\tssum}{{\textstyle\sum}}

\newcommand{\wto}{\rightharpoonup}
\newcommand{\ncdot}{\!\cdot\!}

\newcommand{\gb}{\nabla}
\newcommand{\lec}{\le_c}

\newcommand{\Real}{\mathbb{R}}

\newcommand{\ex}{\mathrm{e}}
\newcommand{\im}{\mathrm{i}}

\newcommand{\tht}{\theta}
\newcommand{\vfi}{\varphi}

\newcommand{\Hper}{H_{\textrm{per}}}

\newcommand{\Dom}{\mathcal{M}}
\newcommand{\Dome}{\mathcal{M}_1}
\newcommand{\Domg}{\mathcal{N}}

\newcommand{\sgb}{\gb^\perp}
\newcommand{\dt}{\;\mathrm{d}t}
\newcommand{\dx}{\;\mathrm{d}\xb}

\newcommand{\kb}{{\boldsymbol{k}}}
\newcommand{\vb}{\boldsymbol{v}}
\newcommand{\xb}{\boldsymbol{x}}

\newcommand{\bmo}{\textsc{bmo}}

\newcommand{\zetah}{\hat\zeta}
\newcommand{\vfih}{\hat\vfi}
\newcommand{\sigmah}{\hat\sigma}
\newcommand{\dex}{\;\mathrm{d}\eta\,\mathrm{d}\xi}

\newcommand{\vrho}{\varrho}

\begin{document}

\title[Mixing rate of 2d perfect fluids]
{A bound on the mixing rate\\ of 2d perfect fluid flows}

\author{D.~Wirosoetisno}
\email{djoko.wirosoetisno@durham.ac.uk}
\urladdr{http://www.maths.dur.ac.uk/\~{}dma0dw}
\address{Mathematical Sciences, Durham University, United Kingdom}


\date{\mydate}

\begin{abstract}
Using the $H^{-1}$ norm as a measure of mixing,
we prove that 2d Euler flows on the torus mix passive scalars
at most exponentially.
The mixing rate is bounded linearly by the BMO norm
of the vorticity (and thus by its $L^\infty$ norm).
We also give an analogous bound on the growth rate of scalar gradients.
\end{abstract}

\maketitle

\section{Introduction}\label{s:intro}

Let $\vb(\xb,t)$ be the velocity field of an incompressible inviscid fluid
flow, governed by the Euler equation
\begin{equation}\label{q:dvdt}
   \dy_t\vb + \vb\ncdot\gb\vb + \gb p = 0
   \qquad\textrm{with}\qquad
   \gb\ncdot\vb = 0.
\end{equation}
For concreteness, we take $(x,y)=\xb\in\Dom:=[0,1]^2$ with periodic boundary
conditions.
With no loss of generality,
we assume that $\vb=(u,v)$ has zero integral over the domain $\Dom$.
For our purpose here we assume that the initial data $\vb(\cdot,0)$ is
smooth, say, $C^3$.
The evolution of a passive scalar $\tht(\xb,t)$ in this flow is then
governed by
\begin{equation}\label{q:dtdt}
   \dy_t\tht + \vb\ncdot\gb\tht = 0.
\end{equation}
As with $\vb$, we assume that the integral of $\tht$ over $\Dom$ vanishes.

Assuming the boundedness of $\vb$ for all $t\ge0$,
one can show from \eqref{q:dtdt} that $|\tht|_{L^p}^{}$ is constant in time
for all $p\in[1,\infty]$.
It is clear that there are many solutions of \eqref{q:dvdt},
e.g., shear flows and (nearly) stable ones,
for which the scalar $\tht$ will undergo little time evolution.
If the solution of \eqref{q:dvdt} is chaotic (as is often the case in
many interesting situations), however,
one expects that $\tht$ will become increasingly ``mixed'' in time;
mathematically, one expects that $\tht(\cdot,t)$ will converge weakly to $0$
as $t\to\infty$ in some suitable space.
Of significant mathematical and physical interest is the rate of this
mixing or convergence.

In \cite{mathew-mezic-petzold:05}, the authors introduced a coarse-grained
``mix norm'' which they showed to be equivalent to the $H^{-1/2}$ norm,
defined formally as
\begin{equation}
   |w|_{H^s}^2 := \tssum_{\kb\ne0}\,|\kb|^{2s}|\hat w_\kb|^2
\end{equation}
where $\hat w_\kb$ are the Fourier coefficients of $w$.
(This is a seminorm in general, but is a proper norm when $w$ has zero
integral over $\Dom$, which will always be the case in this work.)
They argued that this $H^{-1/2}$ norm is a good measure of mixing as it is
weighted towards larger-scale features at the expense of fine details.

Subsequently, \cite{lin-thiffeault-doering:11} proposed the use of the
more convenient $H^{-1}$ norm to measure the mixing rate, having proved that
for any $s<0$ and for any $\tht(\cdot,t)\in L^2$ with zero integral,
the weak convergence $\tht(\cdot,t)\wto0$ in $L^2$ as $t\to\infty$
is equivalent to $|\tht(\cdot,t)|_{H^s}^{}\to0$.

In this work, we follow their use of the $H^{-1}$ norm as our measure of mixing.
Putting $\vfi:=\Delta^{-1}\tht$ and
multiplying \eqref{q:dtdt} by $-\vfi$ in $L^2$, we have
\begin{equation}\label{q:dt1dt}
   \frac12\ddt{\;}|\gb\vfi|_{L^2}^2 - (\vb\ncdot\gb \Delta\vfi,\vfi)_{L^2}^{} = 0.
\end{equation}
Assuming sufficient smoothness, we integrate the nonlinear term by parts
\begin{equation}\label{q:jac}\begin{aligned}
   -(\vb\ncdot\gb\Delta\vfi,\vfi)_{L^2}^{}
	&= \tssum_j\,(\vb\ncdot\gb\dy_j\vfi,\dy_j\vfi)_{L^2}^{}
	+ \tssum_j\,((\dy_j\vb)\ncdot\gb\vfi,\dy_j\vfi)_{L^2}^{}\\
	&= \tssum_j\,((\dy_j\vb)\ncdot\gb\vfi,\dy_j\vfi)_{L^2}^{}.\\
\end{aligned}\end{equation}
Using the standard estimate
\begin{equation}\label{q:naive}
   \bigl|\bigl((\dy_j\vb)\ncdot\gb\vfi,\dy_j\vfi\bigr)_{L^2}^{}\bigr|
	\lec |\gb\vb|_{L^\infty}^{}|\gb\vfi|_{L^2}^2,
\end{equation}
we have
\begin{equation}\label{q:dndt}
   \ddt{\;}|\gb\vfi|_{L^2}^2 \ge -c_1^{}\,|\gb\vb|_{L^\infty}^{}|\gb\vfi|_{L^2}^2.
\end{equation}
Without assuming that $\vb$ is a solution of \eqref{q:dvdt},
the possibility of perfect mixing in a finite time $T$ was not ruled out
in \cite{lin-thiffeault-doering:11}, where lower bounds for perfect mixing
time were given for physically relevant cases;
here perfect mixing is defined by
$|\tht(\cdot,T)|_{H^{-1}}^{}=|\gb\vfi(\cdot,T)|_{L^2}^{}=0$.
They did however point out that if $|\gb\vb(\cdot,t)|_{L^\infty}^{}$ is
uniformly bounded, this will give a bound on the mixing rate.

\medskip
When $\vb$ is the solution of \eqref{q:dvdt}, however, one can use
the argument of \cite{beale-kato-majda:84,kato:86} to rule out perfect
mixing as long as the solution remains regular.
Unlike the results in the following sections,
this argument works both in 2d and 3d.
We fix $s\ge3$ and assume that $\vb(\cdot,0)\in H^s$.
By Sobolev embedding, this implies
$\gb\vb(\cdot,0)\in L^\infty$.
Using the differential inequality
\begin{equation}
   \ddt{\;}|\vb|_{H^s}^2 \le c_2^{}\,|\gb\vb|_{L^\infty}^{}|\vb|_{H^s}^2,
\end{equation}
we have
\begin{equation}\label{q:bdvs}
   |\vb(\cdot,t)|_{H^s}^2 \le |\vb(\cdot,0)|_{H^s}^2 \exp\Bigl(c_2^{}\int_0^t |\gb\vb(\cdot,t')|_{L^\infty}^{} \dt'\Bigr).
\end{equation}
Now \eqref{q:dndt} can also be integrated to give
\begin{equation}\label{q:bdn1}
   |\gb\vfi(\cdot,t)|_{L^2}^2 \ge |\gb\vfi(\cdot,0)|_{L^2}^2 \exp\Bigl(-c_1^{}\int_0^t |\gb\vb(\cdot,t')|_{L^\infty}^{} \dt'\Bigr).
\end{equation}
In the 3d case \cite{beale-kato-majda:84},
one can bound $|\gb\vb|_{L^\infty}^{}$ essentially by
$|\textrm{curl}\,\vb|_{L^\infty}^{}\log|\vb|_{H^3}^{}$ to show that the integral in
\eqref{q:bdvs} is finite for all $t\ge0$ such that
$\textrm{curl}\,\vb\in L^1([0,t];L^\infty)$.
In other words, no perfect mixing is possible before the blow-up time if the
latter is finite.
In the 2d case, it is shown in \cite{kato:86} that $\vb(\cdot,t)\in H^s$
for all $t\ge0$, implying that the integral in \eqref{q:bdn1} is finite
and $|\gb\vfi(\cdot,t)|_{L^2}^{}>0$ for all $t\ge0$.
Therefore, no perfect mixing is possible in finite time when $\vb$ is a
solution of the 2d Euler equation \eqref{q:dvdt}.

\section{A Uniform Bound on the Mixing Rate}\label{s:main}

The inequality \eqref{q:bdn1} does not tell us much about the mixing
{\em rate\/} since we do not have a bound uniform in $t$ for
$|\gb\vb|_{L^\infty}^{}$.
As noted above, \cite{beale-kato-majda:84} bounded $|\gb\vb|_{L^\infty}^{}$
essentially by $|\textrm{curl}\,\vb|_{L^\infty}^{}\log|\vb|_{H^3}^{}$,
which is not bounded uniformly in $t$---in fact, one expects
$|\vb(\cdot,t)|_{H^3}^{}$ to grow without bound for ``generic'' flows.

In this work, we take advantage of the special structure of the Jacobian to
bound \eqref{q:naive} by $|\gb\vb|_\bmo^{}$ instead of $|\gb\vb|_{L^\infty}^{}$
and use a time-uniform bound on $|\gb\vb|_\bmo^{}$ to obtain a bound on the
mixing rate.
To this end,
we introduce the vorticity $\omega:=\dy_xv-\dy_yu$, whose integral over $\Dom$
vanishes by Stokes' theorem, and
the streamfunction $\psi:=\Delta^{-1}\omega$ (defined uniquely by requiring
that its integral over $\Dom$ vanishes), in terms of which
$\vb=\sgb\psi=(-\dy_y\psi,\dy_x\psi)$.
Taking the curl of \eqref{q:dvdt}, we have
\begin{equation}\label{q:dwdt}
   \dy_t\omega + \dy(\psi,\omega) = 0
\end{equation}
where the Jacobian $\dy(\psi,\omega):=\sgb\psi\ncdot\gb\omega=\dy_x\psi\,\dy_y\omega-\dy_x\omega\,\dy_y\psi$.
Using \eqref{q:jac}, we can write \eqref{q:dt1dt} as
\begin{equation}\label{q:dt2dt}
   \frac12\ddt{\;}|\gb\vfi|_{L^2}^2 + \tssum_j\,\bigl(\dy(\dy_j\psi,\vfi),\dy_j\vfi\bigr)_{L^2}^{} = 0.
\end{equation}

For our purpose here, the following definition suffices;
we refer the reader to \cite{stein:ha} for more details.
Let $\Domg\subset\Real^2$.
For $w\in L_{\textrm{loc}}^1(\Domg)$ and any ball $B\subset\Domg$, let
$w_B^{}$ denote the average of $w$ in $B$,
\begin{equation}
   w_B^{} := \frac1{|B|} \int_B w(\xb) \dx.
\end{equation}
We say that $w\in{}$BMO$(\Domg)$, the space of functions of bounded mean
oscillations in $\Domg$, if
\begin{equation}\label{q:bmodef}
   |w|_\bmo^{} := \sup_{B\subset\Domg}\,\frac1{|B|} \int_B |w(\xb)-w_B^{}| \dx < \infty.
\end{equation}
We note that $|\cdot|_\bmo^{}$ as defined here is a seminorm (constants have
BMO-norm zero), but since we only deal with functions of zero average,
\eqref{q:bmodef} defines a proper norm for all relevant quantities.
We also note that \eqref{q:bmodef} implies
\begin{equation}\label{q:bmolinf}
   |w|_\bmo^{} \le |w|_{L^\infty}^{}.
\end{equation}

\medskip
Our main result is the following:

\begin{theorem}\label{t:main}
Let the passive scalar $\tht$ be the solution of \eqref{q:dtdt} and
$\vb$ that of the Euler equation \eqref{q:dvdt} with initial vorticity
$\omega(\cdot,0)\in L^\infty(\Dom)$.
Then $\tht$ satisfies
\begin{equation}\label{q:bdi}
   |\tht(\cdot,t)|_{H^{-1}}^2 \ge |\tht(\cdot,0)|_{H^{-1}}^2 \exp\Bigl(-\lambda\int_0^t |\omega(\cdot,t')|_\bmo^{} \dt'\Bigr)
\end{equation}
for some constant $\lambda$ depending only on the domain $\Dom$.
Moreover, the mixing rate is bounded by the sup-norm of the initial vorticity as
\begin{equation}\label{q:bd0}
   |\tht(\cdot,t)|_{H^{-1}}^2 \ge |\tht(\cdot,0)|_{H^{-1}}^2 \exp\bigl(-t\lambda\,|\omega(\cdot,0)|_{L^\infty}^{}\bigr).
\end{equation}
\end{theorem}

\smallskip\noindent
We note that our result does not address the more difficult issue
of the existence of exponentially mixing solutions of \eqref{q:dvdt},
or whether the bound \eqref{q:bd0} is attained even qualitatively.
See \cite{haynes-jv:05} and references therein for further discussion.

\begin{proof}
Parts of this proof were inspired by \cite{kozono-taniuchi:00}.
As usual, $c$ denotes a generic positive constant whose value may
differ each time the symbol appears.
Assuming for now the following estimate for the Jacobian,
\begin{equation}\label{q:jacbmo}
   \bigl|\dy(\zeta,\vfi)\bigr|_{L^2}^{}
	\le c\,|\gb\zeta|_\bmo^{}|\gb\vfi|_{L^2}^{},
\end{equation}
we take $\zeta=\dy_x\psi$ and $\zeta=\dy_y\psi$ in turn in \eqref{q:dt2dt}
to get [cf.~\eqref{q:dndt}]
\begin{equation}\label{q:dnndt}
   \ddt{\;}|\gb\vfi|_{L^2}^2 \ge -c\,|\gb\vb|_\bmo^{}|\gb\vfi|_{L^2}^2.
\end{equation}
This implies the analogue of \eqref{q:bdn1} with BMO in place of $L^\infty$.
Unlike the $L^\infty$ case where $|\omega|_{L^\infty}^{}$ does not bound
$|\gb\vb|_{L^\infty}^{}$, here we have
\begin{equation}\label{q:vwbmo}
   |\gb\vb|_\bmo^{} \le c\,|\omega|_\bmo^{},
\end{equation}
whose proof will follow shortly.
Using \eqref{q:vwbmo} in \eqref{q:dnndt} and integrating gives us \eqref{q:bdi}.
Using \eqref{q:bmolinf} and the fact that
$|\omega(\cdot,t)|_{L^\infty}^{}=|\omega(\cdot,0)|_{L^\infty}^{}$
gives us \eqref{q:bd0}.

\medskip
To prove \eqref{q:vwbmo}, we start with the identity
(cf.~\cite[p.~59]{stein:sidpf})
\begin{equation}
   \gb\vb = \gb\sgb\Delta^{-1}\omega = (R_x,R_y)(R_y,-R_x)\,\omega
\end{equation}
where $R_x:=\dy_x(-\Delta)^{-1/2}$ and analogously for $R_y$.
It is shown in \cite[p.~138]{stein:ha} that the Riesz transforms
$(R_x,R_y)$ are bounded operators in the Hardy space $\mathcal{H}^1(\Real^2)$,
which is dual to BMO$(\Real^2)$.
Boundedness of $R:=(R_x,R_y)$ in BMO$(\Real^2)$ then follows by duality.
From the definition \eqref{q:bmodef}, it is clear that the norm in
BMO$(\Dom)$ is identical to that for periodic functions in BMO$(\Real^2)$.

\medskip
We now prove the Jacobian estimate \eqref{q:jacbmo}.
For bounded domains in $\Real^2$ with Dirichlet boundary conditions,
although not stated explicitly,
the proof is essentially contained in \cite{kim:09}.
Its extension to periodic boundary conditions is essentially done in
\cite{gtw3:dodu}.
We reproduce these proofs here (with minor modifications) for convenience.

Let $\Dome:=\cup_{\xb\in\Dom} B_1(\xb)$ where $B_1(\xb)\subset\Real^2$
is the ball of unit radius centred at $\xb$.
We first prove \eqref{q:jacbmo} for
$\zeta\in H^1_0(\Dome)\cap\{\gb\zeta\in{}$BMO$(\Dome)\}$
and $\vfi\in H^1_0(\Dome)$.
Denoting the Fourier transform by a h{\^a}t,
we use the following result from \cite[p.~154]{coifman-meyer:opd}.

Let $\sigmah(\cdot,\cdot)\in C^\infty(\Real^2\times\Real^2-\{0\})$ satisfies
\begin{equation}\label{q:cm1}
   \bigl|\dy_\xi^\alpha\dy_\eta^\beta\hat\sigma(\xi,\eta)\bigr| \le C_{\alpha\beta}(|\xi|+|\eta|)^{-|\alpha|-|\beta|}
   \qquad\forall\,\xi,\,\eta\in\Real^2-\{0\}
\end{equation}
for every multi-indices $\alpha$ and $\beta$, and
\begin{equation}\label{q:cm2}
   \sigmah(0,\eta) = 0.
\end{equation}
Then the bilinear operator
\begin{equation}
   \sigma(f,g)(\xb) := \int_{\Real^2\times\Real^2} \ex^{\im\xb\cdot(\xi+\eta)}\,\sigmah(\xi,\eta)\,\hat f(\xi)\,\hat g(\eta) \dex
\end{equation}
is bounded in $L^2(\Dome)$ as
\begin{equation}\label{q:cm3}
   |\sigma(f,g)|_{L^2}^{} \le c\,|f|_\bmo^{}|g|_{L^2}^{}.
\end{equation}

For the Jacobian $\dy(\zeta,\vfi)$, we have
$\sigmah(\xi,\eta)=-\xi^\perp\cdot\eta$, so in order to satisfy
\eqref{q:cm1}--\eqref{q:cm2} we split $\dy(\zeta,\vfi)$ as follows.
Let $\rho\in C^\infty$ be monotone increasing with $\rho(t)=0$ for $t\le1$
and $\rho(t)=1$ for $t\ge2$, and let $\dy(\zeta,\vfi)=J^>+J^<$ where
\begin{equation}\begin{aligned}
   &J^>(\xb) = -\frac{1}{2\pi} \int_{\Real^2\times\Real^2} \ex^{\im\xb\cdot(\xi+\eta)}\,\rho(|\xi|)\,(\xi^\perp\ncdot\eta)\,\zetah(\xi)\,\vfih(\eta) \dex\\
   &J^<(\xb) = -\frac{1}{2\pi} \int_{\Real^2\times\Real^2} \ex^{\im\xb\cdot(\xi+\eta)}\,[1-\rho(|\xi|)]\,(\xi^\perp\ncdot\eta)\,\zetah(\xi)\,\vfih(\eta) \dex.
\end{aligned}\end{equation}
Taking $\sigmah(\xi,\eta)=\rho(|\xi|)$, we have from \eqref{q:cm3}
\begin{equation}
   |J^>|_{L^2}^{} \le c\,|\gb\zeta|_\bmo^{}|\gb\vfi|_{L^2}.
\end{equation}
And taking
\begin{equation}
   \sigmah(\xi,\eta) = \frac{1-\rho(|\xi|)}{1+|\xi|}\,\xi^\perp,
\end{equation}
\eqref{q:cm3} gives us
\begin{equation}
   |J^<|_{L^2}^{} \le c\,|\gb\zeta|_\bmo^{}|\gb\vfi|_{L^2}.
\end{equation}

Having proved \eqref{q:jacbmo} in $\Dome$ with Dirichlet boundary conditions,
we now extend it to the periodic case,
taking $\zeta\in \Hper^1(\Real^2)\cap \{\gb\zeta\in\textrm{BMO}(\Real^2)\}$
and $\vfi\in \Hper^1(\Real^2)$.
Let $\vrho\in C^\infty(\Real^2;[0,1])$ with
$\vrho=1$ in $\Dom$ and $\vrho=0$ in $\Real^2-\Dome$.
By the above, we have
\begin{equation}
   \bigl|\dy(\vrho\zeta,\vrho\vfi)\bigr|_{L^2(\Dome)}^{}
	\le C(\Dom,\Dome)\,|\gb(\vrho\zeta)|_{\bmo(\Dome)}^{}|\gb(\vrho\vfi)|_{L^2(\Dome)}^{}.
\end{equation}
By the periodicity of $\zeta$ and smoothness of $\vrho$, we have
\begin{equation}\label{q:rho1}
   |\gb(\vrho\zeta)|_{\bmo(\Dome)}^{}
	\le C(\Dom,\Dome,\vrho)\,|\gb\zeta|_{\bmo(\Dom)}^{},
\end{equation}
and similarly
\begin{equation}\label{q:rho2}
   |\gb(\vrho\vfi)|_{L^2(\Dome)}^{}
	\le C(\Dom,\Dome,\vrho)\,|\gb\vfi|_{L^2(\Dom)}^{}.
\end{equation}
One then takes the infimum over $\vrho$ in \eqref{q:rho1}--\eqref{q:rho2}
to remove the dependence of the constants on $\vrho$.
The conclusion follows from these and
\begin{equation}
   \bigl|\dy(\zeta,\vfi)\bigr|_{L^2(\Dom)}^{}
   = \bigl|\dy(\vrho\zeta,\vrho\vfi)\bigr|_{L^2(\Dom)}^{}
   \le \bigl|\dy(\vrho\zeta,\vrho\vfi)\bigr|_{L^2(\Dome)}^{}.
\end{equation}
\end{proof}

\section{Further Remarks}

One can obtain an upper bound on the scalar gradient using the same method:
Multiplying \eqref{q:dtdt} by $-\Delta\tht$ and estimating, we find
[cf.~\eqref{q:jac} and \eqref{q:dnndt}]
\begin{equation}\begin{aligned}
   \ddt{\;}|\gb\tht|_{L^2}^2
	&= -2\,\tssum_j\,\bigl((\dy_j\vb)\ncdot\gb\tht,\dy_j\tht)_{L^2}^{}\\
	&\le C(\Dom)\,|\gb\vb|_\bmo^{}|\gb\tht|_{L^2}^2.
\end{aligned}\end{equation}
As before, since $|\gb\vb|_\bmo^{} \le c\,|\omega|_{L^\infty}^{}$ and using
the fact that the latter quantity is time-invariant, we have [cf.~\eqref{q:bd0}]
\begin{equation}
   |\gb\tht(\cdot,t)|_{L^2}^2 \le |\gb\tht(\cdot,0)|_{L^2}^2 \exp\bigl(\, t\lambda(\Dom) |\omega(\cdot,0)|_{L^\infty}^{} \bigr).
\end{equation}

Since the vorticity equation \eqref{q:dwdt}, or equivalently,
\begin{equation}
   \dy_t\omega + \vb\ncdot\gb\omega = 0,
\end{equation}
is of the same form as \eqref{q:dtdt},
these bounds evidently apply equally well to $|\gb\omega(\cdot,t)|_{L^2}^{}$.
For this {\em active\/} scalar, one expects $\gb\omega$ to grow without bound;
see, e.g., \cite{morgulis-shnirelman-yudovich:08} for estimates in H{\"o}lder
spaces.

\bigskip
\begin{ack}
The author thanks Jacques Vanneste and Xiaoming Wang for insightful
discussions leading to this work.
\end{ack}

\bigskip\hbox to\hsize{\qquad\hrulefill\qquad}\medskip

\end{document}